\documentclass[a4paper]{article}
\usepackage{amsmath,amssymb,amsfonts,amsthm} 
\usepackage{graphicx} 

\usepackage{wrapfig} \usepackage{url} \usepackage[usenames]{color}
\def\b1{\mbox{\boldmath $1$}}

\def\qed{\hfill \rule{4pt}{7pt}}

\def\pf{\noindent {\it{Proof.}\hskip 2pt}}

\theoremstyle{plain}
\newtheorem{thm}{Theorem}[section]
\newtheorem{lem}[thm]{Lemma}

\newtheorem{remark}[thm]{Remark}

\newtheorem{problem}[thm]{Problem}

\makeatletter \@addtoreset{equation}{section} \makeatother \makeatletter
\@addtoreset{thm}{section} \makeatother
\allowdisplaybreaks 

\usepackage[numbers,sort&compress]{natbib} 
\usepackage[bookmarksnumbered, bookmarksopen,
colorlinks,citecolor=blue,linkcolor=blue]{hyperref}

\begin{document}
\begin{center}
{{\large\bf Monotonicity of speed for biased random walk on Galton-Watson tree}\footnote{The project is supported partially by CNNSF (No. 11271204).}}
\end{center}
\begin{center}
Song He, Wang Longmin and Xiang Kainan
\vskip 1mm
\footnotesize{School of Mathematical Sciences, LPMC, Nankai University}\\
\footnotesize{Tianjin City, 300071, P. R. China}\\
\footnotesize{Emails: songhe@mail.nankai.edu.cn} (Song)\\
\footnotesize{~~~~~~wanglm@nankai.edu.cn} (Wang)\\
\footnotesize{~~~~~~~~~~~~kainanxiang@nankai.edu.cn} (Xiang)
\end{center}

\begin{abstract}

Ben Arous, Fribergh and Sidoravicius \cite{GAV2014} proved that speed of biased random walk $RW_\lambda$ on a Galton-Watson tree without leaves is
strictly decreasing for $\lambda\leq \frac{m_1}{1160},$ where $m_1$ is minimal degree of the Galton-Watson tree.
And A\"{\i}d\'{e}kon \cite{EA2013} improved this result to $\lambda\leq \frac{1}{2}.$
In this paper, we prove that for the $RW_{\lambda}$ on a Galton-Watson tree without leaves, its speed is strictly decreasing for
$\lambda\in \left[0,\frac{m_1}{1+\sqrt{1-\frac{1}{m_1}}}\right]$ when $m_1\geq 2;$ and we owe the proof to A\"{\i}d\'{e}kon \cite{EA2013}. 
\end{abstract}

\section{Introduction}
In this paper we study biased random walks $RW_{\lambda}$ on Galton-Watson trees. And we focus on the the following question:
Is speed of $RW_\lambda$ monotonic nonincreasing as a function of its bias $\lambda$ when the Galton-Watson tree has no leaves?

Let $\mathbb{T}$ be a Galton-Watson tree with root $e$, and $\nu$ be its offspring distribution random variable
with $m=E(\nu)>1.$ Denote by $(\Omega, P)$ the associated probability space.
Note $\mathbb{T}$ is super-critical and extinction probability $q=P[\mathbb{T}\ \mbox{is finite}]<1.$
Let $\nu(x)$ be the number of children of a vertex $x\in \mathbb{T}.$
For any $x\in\mathbb{T}\setminus\{e\},$ let $x_{\ast}$ be the parent of $x,$ i.e. the neighbor of $x$ lying on geodesic path from $x$ to $e.$
And write $xi, 1\leq i\leq \nu(x)$ as the children of $x.$

Given Galton-Watson tree $\mathbb{T},$ for any $\lambda\geq 0,$ $\lambda$-biased random walk $RW_{\lambda},$ $(X_n)_{n=0}^{\infty},$ is defined as follows.
The transition probability from $x$ to
an adjacent vertex $y$ is
\[
p(x,y)=\left\{\begin{array}{cl}
\frac{1}{\nu(x)} &{\rm if}\ x=e,\\
      \frac{\lambda}{\lambda+\nu(x)} &{\rm if}\ y=x_{\ast},\ x\neq e,\\
     \frac{1}{\lambda+\nu(x)} &{\rm otherwise}.
   \end{array}
\right.
\]
Clearly $(X_n)_{n=0}^{\infty}$ is a reversible Markov chain for $\lambda>0.$
Let $\mathbf{P}_x$ be the quenched probability of $RW_{\lambda}$ starting at $x$ and
$\mathbb{P}_x$ the annealed probability obtained by the semi-direct product
$\mathbb{P}_x=P\times \mathbf{P}_x.$ Denote the respectively associated expectations by $\mathbf{E}_x$ and $\mathbb{E}_x.$
A motivation for introducing $RW_\lambda$ on trees is that this random walk can be used to obtain almost uniform samples from the set of self-avoiding walks of a given length on a lattice (\cite{AA1985}). And for more motivations on biased random walks on graphs, see surveys
\cite{RRY1997} and \cite{GA2014}.

Lyons \cite{RL1990} showed that there is a critical parameter $\lambda_c$ for $RW_{\lambda}$ on a general tree which is just exponential of Hausdorff dimension of the tree boundary, such that $RW_\lambda$ is transient for $\lambda<\lambda_c$ and recurrent for $\lambda>\lambda_c.$
Then in above paper, Lyons proved that for almost every Galton-Watson tree conditioned on non-extinction, $RW_{\lambda}$ is transient for $0\leq \lambda <m.$  And from Lyons \cite{RL1992}, conditionally on non-extinction, $RW_{m}$ is null recurrent and $RW_\lambda$ is positive recurrent when $\lambda>m.$

Let $\vert x\vert$ be the graph distance between $x$ and $e$ for any vertex $x\in \mathbb{T}.$
Note $\vert x\vert$ is also the generation of $x.$ {\it Fix} $X_0=e.$
Speed $\ell_\lambda$ of $RW_\lambda$ is the almost sure limit (if it exists) of $\frac{|X_n|}{n}$ as $n\rightarrow \infty.$
In this paper, dependence of $\ell_\lambda$ with respect to environment will often be omitted.

Transient $RW_\lambda$ can have zero speed when too much time is spent at leaves.
In \cite{RRY1996}, Lyons, Pemantle and Peres proved that, conditionally on nonextinction, $\ell_\lambda$ exists almost surely, and
$$\ell_{\lambda}\ \mbox{is determinist and positive iff}\ \lambda\in \left(E\left(\nu q^{\nu-1}\right), m\right).$$
From Lyons, Pemantle and Peres \cite{RRY1995}, $\ell_1=E\left(\frac{\nu-1}{\nu+1}\right).$
And in \cite{EA2014},  A\"{\i}d\'{e}kon gave an expression of $\ell_\lambda$ specified in (\ref{eq1}) though an artificial parent to $e$ was added therein.
For related results, refer to Gantert et al. \cite{NSSM2012}.

Lyons, Pemantle and Peres \cite{RRY1996} (see also \cite{RRY1997}) raised the following problem, which was called Lyons-Pemantle-Peres monotonicity problem
by \cite{GAV2014}.

\begin{problem}\label{prob1}
Assume $P(\nu=0)=0,$ namely Galton-Watson tree $\mathbb{T}$ has no leaf.
Is the speed $\ell_\lambda$ of $RW_{\lambda}$ on $\mathbb{T}$ monotonic nonincreasing in
$\lambda\in [0, m)?$
\end{problem}
It was conjectured in \cite{RRY1996} and \cite{RRY1997} that Problem \ref{prob1} should have a positive answer.
Obviously, the answer is positive when $P(\nu=k)=1$ for some $k.$
Whatever, it seems that the speed is nonincreasing for any tree.
But this is wrong. For instance, on binary tree with pipes, which is a multi-type Galton-Watson tree,
the speed is $\frac{(2-\lambda)(\lambda-1)}{\lambda^2+3\lambda-2}$ for $1\leq \lambda\leq 2.$
And also for any $0<\lambda_1<\lambda_2,$ by the repeated filtering method, one can produce a tree such that the speed of $RW_{\lambda_1}$
is less than that of $RW_{\lambda_2}.$ Refer to \cite{RRY1997} for these facts. Notice the just mentioned examples are not
Galton-Watson trees and show the complexity of Problem \ref{prob1}.
Therefore, if the monotonicity of $\ell_\lambda$ holds, then it will be a very fundamental special property of Galton-Watson trees.

For Galton-Watson trees without leaves, the Lyons-Pemantle-Peres monotonicity problem was answered positively for
$\lambda\leq \frac{m_1}{1160}$ by Ben Arous, Fribergh and Sidoravicius \cite{GAV2014},
where
$$m_1=\min\{k\geq 1:\ P[\nu=k]>0\}$$
is minimal degree of the Galton-Watson tree.
And A\"{\i}d\'{e}kon \cite{EA2013} improved the just mentioned result to $\lambda\leq \frac{1}{2}$ by a
completely different approach. In \cite{GYSO2013}, Ben Arous, Hu, Olla and Zeitouni obtained the Einstein
relation for $RW_{\lambda}$ on Galton-Watson trees, which implies
Problem \ref{prob1} holds in a neighborhood of $m.$ These very slow progresses show Problem \ref{prob1} is rather difficult.
For more information on $RW_\lambda$ on $\mathbb{T},$
see \cite{RRY1997} and \cite{GA2014} and references therein. And for monotonicity of speed of biased random walk on groups,
see \cite{Song2015}.

Now our main result is stated as follows.
\begin{thm}\label{thm1}
The speed $\ell_\lambda$ of $RW_{\lambda}$ on Galton-Watson tree $\mathbb{T}$ without leaves
is strictly decreasing in $\lambda\in \left[0,\frac{m_1}{1+\sqrt{1-\frac{1}{m_1}}}\right]$ when $m_1\geq 2.$
\end{thm}
\vskip 2mm
\begin{remark}
It is interesting to provide a nontrivial example to confirm Problem \ref{prob1}. This is an embarrassing problem
in a certain sense. The upper bound of $\lambda$ in Theorem \ref{thm1} is not optimal. And we do not think we can
improve it to $m_1,$ not to mention $m.$ To answer Problem \ref{prob1} affirmatively, one potential approach is firstly to improve
Lemma \ref{lem2} as follows: almost surely, for some positive random variable $c_{\mathbb{T}_\ast},$
$\beta_{\mathbb{T}_\ast}(e,\lambda)
   \geq c_{\mathbb{T}_\ast}(m-\lambda),\ \forall\lambda<m;$
and then to refine the method of present paper and \cite{EA2013}.
\end{remark}

\section{Proof of Theorem \ref{thm1}}\label{sec2}
Inspired by \cite{EA2013}, based on some new observations, we prove Theorem \ref{thm1}.

Let $\mathbb{T}_{\ast}$ be the tree obtained from $\mathbb{T}$ by adding an artificial parent $e_{\ast}$ to the root $e.$
For any vertex $x \in \mathbb{T}_{\ast}$, let
\[\tau_x=\min\{n\geq 0, X_n=x\},\]
where $\min \emptyset=\infty,$ and $(X_n)_{n=0}^\infty$ is a $\lambda$-biased random walk on $\mathbb{T}_{\ast}.$
And for $x \neq e_{\ast},$ let
\[\beta(x):=\beta(x,\lambda)=\mathbf{P}_x(\tau_{x_{\ast}}=\infty)\]
be the quenched probability of never reaching the parent $x_{\ast}$ of $x$ when starting from $x.$
Since $\mathbb{T}$ has no leaf and $\lambda<m,$ we have $\beta(x)>0$ due to transience.
Let $(\beta_i)_{i\geq 0}$ be generic $i.i.d.$ random variables distributed under $\mathbb{P}$ as $\beta(e),$ and
independent of $\nu.$

In \cite{EA2014}, A\"{\i}d\'{e}kon gave the following expression of $\ell_{\lambda}:$
\begin{eqnarray}\label{eq1}
\ell_{\lambda}=\mathbb{E}\left(\frac{(\nu-\lambda)\beta_0}{\lambda-1+\sum_{i=0}^\nu\beta_i}\right)/
    \mathbb{E}\left(\frac{(\nu+\lambda)\beta_0}{\lambda-1+\sum_{i=0}^\nu\beta_i}\right),
    \ \lambda<m.
\end{eqnarray}
Notice (\ref{eq1}) holds trivially for $\lambda=0.$ Here we point out that $RW_{\lambda}$ on $\mathbb{T}_{\ast}$ and $RW_\lambda$ on $\mathbb{T}$ has a slight difference, but due to $\lambda<m$ and transience, these two biased random walks have the same speed when starting at $e.$
{\it Indeed}, we have the following

\begin{lem}\label{lem0}
For $\lambda<m,$ $RW_{\lambda}$ on $\mathbb{T}_{\ast}$ and $RW_\lambda$ on $\mathbb{T}$ have the same speed when starting at $e.$
\end{lem}
\pf For $RW_\lambda$ $(X_n)_{n=0}^\infty$ on $\mathbb{T}_{\ast}$ with $X_0=e,$ define
\begin{eqnarray*}
&&\tau_0=0,\ \sigma_0=\inf\{n\geq\tau_0:\ X_n\not\in\{e,e_{\ast}\}\};\\
&&\tau_1=\inf\{n\geq \sigma_0:\ X_n=e\},\
   \sigma_1=\inf\{n\geq\tau_1:\ X_n\not\in\{e,e_{\ast}\}\}\ \mbox{when}\ \tau_1<\infty;\\
&&\mbox{and for any}\ i\geq 1,\\
&&\tau_{i+1}=\inf\{n\geq \sigma_i:\ X_n=e\},\
   \sigma_{i+1}=\inf\{n\geq\tau_{i+1}:\ X_n\not\in\{e,e_{\ast}\}\}\ \mbox{when}\ \tau_{i+1}<\infty.
\end{eqnarray*}
Since $RW_\lambda$ $(X_n)_{n=0}^\infty$ is transient, there is a unique $i_\ast$ such that $\tau_{i_\ast}<\infty$ and $\tau_{i_\ast+1}=\infty.$
Define a random walk $(Y_n)_{n=0}^{\infty}$ as follows:
\begin{eqnarray*}
(Y_n)_{n=0}^{\infty}=\left(X_{\tau_0},\underbrace{X_{\sigma_0},\cdots,X_{\tau_1}},\underbrace{X_{\sigma_1},\cdots,X_{\tau_2}},\cdots,
    \underbrace{X_{\sigma_{i_\ast-1}},\cdots,X_{\tau_{i_\ast}}},X_{\sigma_{i_\ast}},X_{\sigma_{i_\ast}+1},\cdots\right).
\end{eqnarray*}
Then it is easy to see that $(Y_n)_{n=0}^{\infty}$ is just an $RW_\lambda$ on $\mathbb{T}$ starting at $e.$
It is known that almost surely, both $\lim\limits_{n\rightarrow\infty}\frac{\vert Y_n\vert}{n}$ and $\lim\limits_{n\rightarrow\infty}\frac{\vert X_n\vert}{n}$
exist and are deterministic. By our construction, there exists a random function $s(\cdot)$ on nonnegative integers such that almost surely,
$$Y_n=X_{s(n)},\ n\geq 1,\ \mbox{and}\ \lim\limits_{k\rightarrow\infty}\frac{s(k)}{k}=1.$$
Therefore, almost surely,
\[\lim_{n\rightarrow\infty}\frac{\vert Y_n\vert}{n}=\lim_{n\rightarrow\infty}\frac{\vert X_{s(n)}\vert}{n}
   =\lim_{n\rightarrow\infty}\frac{\vert X_{s(n)}\vert}{s(n)}=\lim_{n\rightarrow\infty}\frac{\vert X_{n}\vert}{n}.\]
This implies the lemma.\qed
\vskip 2mm

For any $n\geq 1,$ let $\beta_n(x):=\beta_n(x,\lambda)$ be the probability to hit level $n$ before $x_{\ast}$ when $\vert x\vert\leq n.$
Recall for vertex $x,$ $xi$ is its $i$-th child and $\nu(x)$ is the number of its children.
Then
$\beta_n(x)=1$ if $|x|=n;$ and for $|x|<n,$
\begin{eqnarray} \label{eq2}
\beta_n(x)=\frac{\sum_{i=1}^{v(x)}\beta_n(xi)}{\lambda+\sum_{i=1}^{v(x)}\beta_n(xi)}.
\end{eqnarray}
Clearly, $\beta_n(e)\downarrow\beta(e)$ as $n\uparrow\infty,\ a.s.,$ and each $\beta_n(x)$ has a continuous derivative in $\lambda$ when $|x|\leq n.$
Put
\begin{eqnarray*}
&&A_n(x)=\frac{\lambda}{\left(\lambda+\sum_{i=1}^{v(x)}\beta_n(xi)\right)^2},\\
&&B_n(x)=\frac{\sum_{i=1}^{v(x)}\beta_n(xi)}{\left(\lambda+\sum_{i=1}^{v(x)}\beta_n(xi)\right)^2}.
\end{eqnarray*}

To continue, we need the following Lemmas \ref{lem1} and \ref{lem2}. For any natural number $d,$ let $\mathbb{T}_{d+1}$
be the $d+1$-regular tree. Define the following generating function
\[U(x,y|z)=\sum_{n=0}^{\infty}\mathbf{P}_x\left(\tau_y^{+}=n\right)z^n,\ x,y\in\mathbb{T}_{d+1},\ -1\leq z\leq 1,\]
where $\tau_y^{+}$ is the first positive time hitting $y$ and $\mathbf{P}_x$ is the law of $RW_\lambda$ starting at $x$ on $\mathbb{T}_{d+1}$ with
a fixed root $e.$
Clearly $\mathbf{P}_x\left(\tau_x^{+}=0\right)=0.$

\begin{lem}\label{lem1}
For any $\lambda\geq 0$ and any vertex $x\in \mathbb{T}_{d+1}\setminus\{e\}$ with parent $y,$
\[\mathbf{P}_x\left(\tau_y^{+}=\infty\right)=1-U(x,y|1)=1-\frac{\lambda\wedge d}{d}.\]
\end{lem}
\pf Obviously $U(x,y|z)$ is absolutely convergent for $|z|\leq 1.$ And for $0\leq z\leq 1,$ we have a probability interpretation:
$U(x,y|z)$ is the probability of ever visiting $y$ in the random walk where it dies out at each
step with probability $1-z.$

Consider the $\lambda$-biased random walk on $\mathbb{T}_{d+1}.$
Since $y$ is the parent of $x,$ by taking one step on $\mathbb{T}_{d+1}$ starting at $x,$ we
can see that either with probability $\frac{\lambda}{\lambda+d}$ the random walk hits $y,$
or the random walk moves to the children of $x$ with probability $\frac{d}{\lambda+d}.$
Notice that in the second case, in order to return $y,$ the random walk must return firstly to $x$ and
then hit $y.$ So by the symmetry of $\mathbb{T}_{d+1},$
\[U(x,y|z)=\frac{\lambda}{\lambda+d}z+\frac{d}{\lambda+d}zU(x_1,x|z)U(x,y|z).\]
Here $x_1$ is a child of $x.$ Notice that $\tau_x^{+}$ under $\mathbf{P}_{x_1}$ has the same law as $\tau_y^{+}$ under $\mathbf{P}_x.$
We have that
$$U(x_1,x|z)=U(x,y|z).$$
Therefore,
\[U(x,y|z)=\frac{\lambda}{\lambda+d}z+\frac{d}{\lambda+d}zU(x,y|z)^2;\]
which implies that
\[U(x,y|z)=\frac{(\lambda+d)\pm \sqrt{(\lambda+d)^2-4d \lambda z^2}}{2dz}.\]
Due to $U(x,y|z)$ is continuous for $|z|<1,$ we have
\[U(x,y|z)=\frac{(\lambda+d)- \sqrt{(\lambda+d)^2-4d \lambda z^2}}{2dz}.\]
And further
\[U(x,y|1)=\frac{\lambda+d-\vert d-\lambda\vert}{2d}=\frac{\lambda\wedge d}{d}.\]
By the definition of $\tau_y^{+},$ we obtain
\[\mathbf{P}_x\left(\tau_y^{+}=\infty\right)=1-U(x,y|1)=1-\frac{\lambda\wedge d}{d}.\] \qed
\vskip 2mm
Let us interpret $\beta(e)$ in the framework of electric networks. Given any weighted graph (in another word, electric network)
$G=(V(G),E(G))$ with nonnegative edge weight function $c.$ Note weights are called conductances. Suppose $a\in V(G)$ and $Z\subseteq  V(G).$ Write
\[\mathbf{P}^{G,c}_a(a\rightarrow Z)=\mathbf{P}^{G,c}_a\left(\tau_Z< \tau_a^+\right),\]
where $\tau_Z=\inf\{n\geq 0:\ X_n\in Z\},\ \tau_a^+=\inf\{n\geq 1:\ X_n=a\},$ $(X_n)_{n\geq 0}$ is the random walk associated with electric network
$G,$ and $\mathbf{P}^{G,c}_a$ is the law of $(X_n)_{n\geq 0}$ starting at $a.$ Let
$$\pi(x)=\sum\limits_{y\in V(G):\ y\sim x}c(\{x,y\}),\ \forall x\in V(G),\ \mbox{where}\ y\sim x\ \mbox{means}\ y\ \mbox{is adjacent to}\ x.$$
Then $\pi(\cdot)$ is a stationary measure of $(X_n)_{n\geq 0}.$ Call
$$\mathcal{C}_G(a\leftrightarrow Z):=\mathcal{C}_{G,c}(a\leftrightarrow Z)=\pi(a)\mathbf{P}^{G,c}_a(a\rightarrow Z)$$
effective conductance between $a$ and $Z.$ Use $\mathbf{P}^{G,c}_{a}(a\rightarrow \infty)$ to denote the probability of
$(X_n)_{n\geq 0}$ never returning to $a$ when $X_0=a.$ Then call
$$\mathcal{C}_G(a\leftrightarrow \infty):=\mathcal{C}_{G,c}(a\leftrightarrow \infty)=\pi(a)\mathbf{P}^{G,c}_a(a\rightarrow \infty)$$
effective conductance from $a$ to $\infty$ in $G.$

To emphasize on $\mathbb{T}_\ast,$ denote $\beta(e)=\beta(e,\lambda)$ by $\beta_{\mathbb{T}_\ast}(e,\lambda).$
When $\lambda>0,$ on $\mathbb{T}_\ast,$ endow any edge $\{x,y\}$ with $x,y\not=e_\ast$ with a weight $\lambda^{-\vert x\vert\wedge\vert y\vert-1},$
and edge $\{e_\ast,e\}$ with a weight 1; and denote this weight function by $c_0.$
Then for $\lambda>0,$ the $RW_\lambda$ on $\mathbb{T}_\ast$ is the random walk associated with weighted graph
(electric network) $\mathbb{T}_\ast;$ and
$$\beta_{\mathbb{T}_\ast}(e,\lambda)=\mathbf{P}^{\mathbb{T}_\ast,c_0}_{e_\ast}(e_\ast\rightarrow\infty)
  =\mathcal{C}_{\mathbb{T}_\ast,c_0}(e_\ast\leftrightarrow\infty).$$

\begin{lem}\label{lem2}
Assume Galton-Watson tree $\mathbb{T}$ has no leaf. Then almost surely
\[1-\frac{\lambda\wedge m_1}{m_1} \leq\beta_{\mathbb{T}_\ast}(e,\lambda)\leq 1-\frac{\lambda}{m_2},\ \lambda\in [0,m),\]
where $m_2=\sup\{k\geq 1:\ P[\nu=k]>0\}.$
\end{lem}
\pf When $m_2=\infty,$ $\beta(e)\leq 1-\frac{\lambda}{m_2}$ holds trivially. Clearly the lemma is true for $\lambda=0.$
So we assume $m_2<\infty$ (namely $\nu$ takes finitely many values) and $0<\lambda<m.$
Through a natural way, we can embed an $m_1$-ary tree $\mathbb{H}^1$ into $\mathbb{T}$
and also embed $\mathbb{T}$ into an $m_2$-ary tree $\mathbb{H}^2$ such that roots of $\mathbb{H}^1$ and $\mathbb{H}^2$
are root $e$ of $\mathbb{T}.$
Similarly to $\mathbb{T}_{\ast},$ let each $\mathbb{H}^i_{\ast}$ be obtained from $\mathbb{H}^i$ by adding the artificial parent $e_{\ast}$ of $e$
to $\mathbb{H}^i.$

Like electric network $(\mathbb{T}_\ast,c_0),$ we endow a weight function $c_i$ to each $\mathbb{H}^i_\ast.$
And view $c_0$ and $c_1$ as functions on the set of edges of $\mathbb{H}^2_\ast$ by letting that
$$c_0(\{x,y\})=0\ (\mbox{resp.}\ c_1(\{x,y\})=0)$$
when $\{x,y\}$ is not an edge of $\mathbb{T}_\ast$ (resp. $\mathbb{H}^1_\ast$).
Then $c_1(\{x,y\})\leq c_0(\{x,y\})\leq c_2(\{x,y\})$ for any edge $\{x,y\}$ of $\mathbb{H}^2_\ast.$

Notice that
\begin{eqnarray*}
&&\beta_{\mathbb{T}_\ast}(e,\lambda)=\mathbf{P}^{\mathbb{T}_{\ast},c_0}_{e_\ast}(e_{\ast}\rightarrow
     \infty)=\mathbf{P}^{\mathbb{H}^2_{\ast},c_0}_{e_\ast}(e_{\ast}\rightarrow \infty)=
     \mathcal{C}_{\mathbb{H}^2_\ast,c_0}(e_\ast\leftrightarrow\infty),\\
&&\beta_{\mathbb{H}^i_\ast}(e,\lambda)=\mathbf{P}^{\mathbb{H}^i_\ast,c_i}_{e_\ast}(e_{\ast}\rightarrow \infty)=
    \mathbf{P}^{\mathbb{H}^2_{\ast},c_i}_{e_\ast}(e_{\ast}\rightarrow \infty)=
     \mathcal{C}_{\mathbb{H}^2_\ast,c_i}(e_\ast\leftrightarrow\infty),\ i=1,2.
\end{eqnarray*}

Recall Rayleigh's monotonicity principle from \cite{RLYP2014} Section 2.4:
Let $G$ be an infinite connected graph with two nonnegative edge weight functions $c$ and $c^\prime$
such that $c\leq c^\prime$ everywhere. Then for any vertex $a$ of $G,$
$\mathcal{C}_{G,c}(a\leftrightarrow \infty)\leq \mathcal{C}_{G,c'}(a\leftrightarrow \infty).$

Therefore, we have that
$$\mathcal{C}_{\mathbb{H}^2_\ast,c_1}(e_\ast\leftrightarrow\infty)\leq \mathcal{C}_{\mathbb{H}^2_\ast,c_0}(e_\ast\leftrightarrow\infty)
   \leq \mathcal{C}_{\mathbb{H}^2_\ast,c_2}(e_\ast\leftrightarrow\infty).$$
Namely
$$\beta_{\mathbb{H}^1_\ast}(e,\lambda)\leq \beta_{\mathbb{T}_\ast}(e,\lambda)\leq \beta_{\mathbb{H}^2_\ast}(e,\lambda).$$
Hence by Lemma \ref{lem1}, we obtain that
\begin{eqnarray*}
1-\frac{\lambda\wedge m_1}{m_1}\leq\beta_{\mathbb{T}_\ast}(e,\lambda)\leq 1-\frac{\lambda}{m_2}.
\end{eqnarray*}
The lemma holds.\qed
\vskip 2mm

Now we are in the position to prove the following lemma on derivative of $\beta(e,\lambda):$

\begin{lem}\label{lem3}
For Galton-Watson tree $\mathbb{T}$ without leaves, almost surely, $\beta(e)=\beta(e,\lambda)$ has continuous derivative $\beta^\prime(e)=\beta^\prime(e,\lambda)$ in $\lambda\in [0,m_1),$ and
\begin{eqnarray}
\label{eq3} 0<-\beta^{\prime}(e,\lambda)\leq \frac{\beta(e,\lambda)}{m_1-\lambda},\ \lambda\in [0,m_1).
\end{eqnarray}

\end{lem}
\pf
Derivating (\ref{eq1}) in $\lambda<m$ yields that
\[-\beta_{n}^{\prime}(x,\lambda)=A_n(x)\sum_{i=1}^{\nu(x)}-\beta^{\prime}_n(xi,\lambda)+B_n(x),\]
where $\beta^{\prime}_n(x,\lambda)$ is the derivative in $\lambda.$ Then
\begin{eqnarray}\label{eq4}
-\beta_{n}^{\prime}(e,\lambda)=\sum_{k=0}^{n-1}\sum_{|x|=k}B_n(x)\prod_{i=0}^{k-1}A_n(x_i),\ \lambda<m,
\end{eqnarray}
where $x_i$ is the ancestor at generation $i$ of $x.$
And for any $k\in [0,n-1]$,
\[\beta_n(e,\lambda)=\sum_{|x|=k}\beta_n(x,\lambda)\prod_{i=0}^{k-1}\frac{1}{\lambda+\sum_{i=1}^{\nu(x_i)}\beta_n(x_ij,\lambda)},\ \lambda<m.\]
Here $x_ij$ is the $j$-th child of the ancestor $x_i$.

Clearly, $\beta_n(x,\lambda)$ is nonincreasing in $n.$
By Lemma \ref{lem2}, we have that
\[\lambda+\sum_{i=1}^{\nu(x)}\beta_n(xi,\lambda)\geq m_1,\ \lambda<m.\]
Hence, for $\lambda<m,$
\begin{eqnarray}\label{eq5}
A_n(x)\leq\frac{1}{m_1}\frac{\lambda}{\lambda+\sum_{i=1}^{\nu(x)}\beta_n(xi,\lambda)},\
B_n(x)\leq \frac{1}{m_1} \beta_n(x,\lambda).
\end{eqnarray}
And for $\lambda<m,$
\begin{eqnarray}\label{eq6}
\sum_{|x|=k}B_n(x)\prod_{i=0}^{k-1}A_n(x_i)\leq \frac{1}{m_1^{k+1}}\sum_{|x|=k}\beta_n(x,\lambda)\prod_{i=0}^{k-1}\frac{\lambda}{\lambda+\sum_{i=1}^{\nu(x_i)}\beta_n(x_ij,\lambda)}
  =\frac{\lambda^k}{m_1^k}\frac{1}{m_1}\beta_n(e,\lambda).
\end{eqnarray}
By (\ref{eq4}) and (\ref{eq6}), almost surely,
\begin{eqnarray}\label{eq7}
0\leq -\beta_n^\prime(e,\lambda)\leq \frac{\beta_n(e,\lambda)}{m_1-\lambda}\leq\frac{1}{m_1-\lambda},\ \lambda<m_1.
\end{eqnarray}

Given any small enough $\epsilon>0.$ From (\ref{eq7}), we see that almost surely,
$$\{(\beta_n(e,\lambda):\ \lambda\in [0,m_1-\epsilon])\}_{n\geq 1},$$
as a sequence of functions on $[0,m_1-\epsilon],$ is equi-continuous. Combining with $\beta_n(e,\lambda)\downarrow\beta(e,\lambda)$ as $n\uparrow\infty$
for all $\lambda\in [0,m),\ a.s.,$  by the Ascoli-Arzel\`{a} theorem, $\left\{\beta_n(e,\lambda):\ \lambda\in [0,m_1-\epsilon]\right\}_{n\geq 1}$ converges
uniformly to $(\beta(e,\lambda):\ \lambda\in [0,m_1-\epsilon])$ almost surely.

Note for any vertex $x\in\mathbb{T},$ $\left(\left\{(\beta_n(x,\lambda):\ 0\leq \lambda<m)\right\}_{n\geq 1},\ (\beta(x,\lambda):\ 0\leq\lambda<m)\right)$ has the same distribution as $\left(\left\{(\beta_n(e,\lambda):\ 0\leq \lambda<m)\right\}_{n\geq 1},\ (\beta(e,\lambda):\ 0\leq\lambda<m)\right).$ We obtain that almost surely, for any vertex $x\in\mathbb{T},$
$\left\{\beta_n(x,\lambda):\ \lambda\in [0,m_1-\epsilon]\right\}_{n\geq 1}$ converges
uniformly to $(\beta(x,\lambda):\ \lambda\in [0,m_1-\epsilon]).$
Hence by the definitions of $A_n(x)$ and $B_n(x),$ we have that almost surely, for any vertex $x,$ $A_n(x)$ and $B_n(x)$ converge uniformly in
$\lambda\in [0,m_1-\epsilon]$ to some continuous functions $A(x)$ and $B(x)$ respectively.

Notice (\ref{eq4}) and (\ref{eq6}). By the dominated convergence theorem, we see that almost surely,
$\left(\beta^{\prime}_n(x,\lambda):\ \lambda\in [0,m_1-\epsilon]\right)$
converges uniformly to some continuous function $(F_{\lambda}:\ \lambda\in [0,m_1-\epsilon]).$
And by the dominated convergence theorem again, almost surely,
$\int_0^{\lambda}\beta_{n}^{\prime}(e,s)\ {\rm d}s$ converges to $\int_0^{\lambda}F_s\ {\rm d}s$ which is also $\beta(e,\lambda)-1$
for all $\lambda\leq m_1-\epsilon.$

Since $\epsilon$ is arbitrary, we obtain that almost surely, $\beta(e,\lambda)$ is differentiable in $\lambda\in [0,m_1).$
And further, almost surely,
$$0\leq-\beta^\prime(e,\lambda)\leq\frac{\beta(e)}{m_1-\lambda},\ \lambda\in [0,m_1).$$
By checking (2.4) and definitions of $A_n(x)$ and $B_n(x),$ when taking limits, we indeed have almost surely,
$$0<-\beta^\prime (e,\lambda)\leq\frac{\beta(e)}{m_1-\lambda},\ \lambda\in [0,m_1).$$
\qed
\vskip 2mm

By symmetry, we have that
\[\ell_{\lambda}=\mathbb{E}\left(\frac{\nu-\lambda}{\nu+1}\frac{\sum_{i=0}^\nu\beta_i}{\lambda-1+\sum_{i=0}^\nu\beta_i}\right)/\mathbb{E}\left(\frac{\nu+\lambda}{\nu+1}\frac{\sum_{i=0}^\nu\beta_i}{\lambda-1+\sum_{i=0}^\nu\beta_i}\right).\]
By Lemma \ref{lem3}, each $\beta_i$ has derivative in $\lambda\in [0,m_1),$ and so does $\ell_\lambda.$
Write each $\beta^{\prime}_i$ and $\ell_{\lambda}^\prime$ for the derivatives in $\lambda\in [0,m_1)$ of $\beta_i$ and $\ell_\lambda$ respectively.
Then by a straightforward calculus \cite{EA2013}, for $\lambda\in [0,m_1),$ $\ell_\lambda^{\prime}<0$ is equivalent with
\begin{eqnarray}\label{eq8}
\nonumber
&&\mathbb{E}\left(\frac{\nu}{\nu+1}\frac{\sum_{i=0}^\nu\beta_i}{\lambda-1+\sum_{i=0}^\nu\beta_i}\right)\mathbb{E}\left(\frac{1}{\nu+1}\frac{\sum_{i=0}^\nu(\beta_i+(1-\lambda)\beta_i^{\prime})}{(\lambda-1+\sum_{i=0}^\nu\beta_i)^2}\right)\\
\nonumber
  &&\ \ \ \ - \mathbb{E}\left(\frac{1}{\nu+1}\frac{\sum_{i=0}^\nu\beta_i}{\lambda-1+\sum_{i=0}^\nu\beta_i}\right)\mathbb{E}\left(\frac{\nu}{\nu+1}\frac{\sum_{i=0}^\nu(\beta_i+(1-\lambda)\beta_i^{\prime})}{(\lambda-1+\sum_{i=0}^\nu\beta_i)^2}\right)\\
&&<\frac{1}{\lambda}\mathbb{E}\left(\frac{\nu}{\nu+1}\frac{\sum_{i=0}^\nu\beta_i}
{\lambda-1+\sum_{i=0}^\nu\beta_i}\right)\mathbb{E}\left(\frac{1}{\nu+1}
\frac{\sum_{i=0}^\nu\beta_i}{\lambda-1+\sum_{i=0}^\nu\beta_i}\right).
\end{eqnarray}
\vskip 2mm
\begin{lem}\label{lem4}
For Galton-Watson tree $\mathbb{T}$ without leaves, (\ref{eq8}) is true for $\lambda\in \left[0,\frac{m_1}{1+\sqrt{1-\frac{1}{m_1}}}\right]$
when $m_1\geq 2.$
\end{lem}
\pf Note that by Lemma \ref{lem3},
$$0<-\beta_i^\prime(\lambda)\leq\frac{\beta_i(\lambda)}{m_1-\lambda},\ \lambda<m_1,\ i\geq 0.$$
By Lemma \ref{lem2},
\begin{eqnarray}\label{eq9}
\lambda-1+\sum_{i=0}^{\nu}\beta_i\geq \lambda-1+(1-\lambda/m_1)\times (m_1+1)=m_1-\frac{\lambda}{m_1}>m_1-1,\ \lambda<m_1.
\end{eqnarray}
Then for any $\lambda<1,$
\begin{eqnarray*}
&&0\leq \beta_i(\lambda)+(1-\lambda)\beta_i^\prime(\lambda)<\beta_i(\lambda),\\
&&\mathbb{E}\left(\frac{1}{\nu+1}\frac{\sum_{i=0}^\nu\beta_i}{\lambda-1+\sum_{i=0}^\nu\beta_i}\right)
  \mathbb{E}\left(\frac{\nu}{\nu+1}\frac{\sum_{i=0}^\nu(\beta_i+(1-\lambda)\beta_i^{\prime})}
  {(\lambda-1+\sum_{i=0}^\nu\beta_i)^2}\right)\geq 0.
\end{eqnarray*}
And further when $\lambda<1,$
\begin{eqnarray*}
&&\mathbb{E}\left(\frac{\nu}{\nu+1}\frac{\sum_{i=0}^\nu\beta_i}{\lambda-1+\sum_{i=0}^\nu\beta_i}\right)\mathbb{E}
            \left(\frac{1}{\nu+1}\frac{\sum_{i=0}^\nu(\beta_i+(1-\lambda)\beta^{\prime}_i)}{(\lambda-1+\sum_{i=0}^\nu\beta_i)^2}\right)\\
 &&\ \ \ \ <\frac{1}{m_1-\frac{\lambda}{m_1}}
              \mathbb{E}\left(\frac{\nu}{\nu+1}\frac{\sum_{i=0}^\nu\beta_i}{\lambda-1+\sum_{i=0}^\nu\beta_i}\right)
              \mathbb{E}\left(\frac{1}{\nu+1}\frac{\sum_{i=0}^\nu\beta_i}{\lambda-1+\sum_{i=0}^\nu\beta_i}\right).
\end{eqnarray*}
Since $\lambda<1$ and $m_1\geq 2,$ $m_1-\frac{\lambda}{m_1}>\lambda,$ we obtain that when $\lambda<1,$
\begin{eqnarray*}
&&\mathbb{E}\left(\frac{\nu}{\nu+1}\frac{\sum_{i=0}^\nu\beta_i}{\lambda-1+\sum_{i=0}^\nu\beta_i}\right)\mathbb{E}\left(\frac{1}{\nu+1}\frac{\sum_{i=0}^\nu(\beta_i+(1-\lambda)\beta_i^{\prime})}{(\lambda-1+\sum_{i=0}^\nu\beta_i)^2}\right)\\
&&\ \ \ \ -
    \mathbb{E}\left(\frac{1}{\nu+1}\frac{\sum_{i=0}^\nu\beta_i}{\lambda-1+\sum_{i=0}^\nu\beta_i}\right)\mathbb{E}\left(\frac{\nu}{\nu+1}\frac{\sum_{i=0}^\nu(\beta_i+(1-\lambda)\beta_i^{\prime})}{(\lambda-1+\sum_{i=0}^\nu\beta_i)^2}\right)\\
&&<\frac{1}{\lambda}\mathbb{E}\left(\frac{\nu}{\nu+1}\frac{\sum_{i=0}^\nu\beta_i}
   {\lambda-1+\sum_{i=0}^\nu\beta_i}\right)\mathbb{E}\left(\frac{1}{\nu+1}
    \frac{\sum_{i=0}^\nu\beta_i}{\lambda-1+\sum_{i=0}^\nu\beta_i}\right);
\end{eqnarray*}
namely (\ref{eq8}) holds.

When $\lambda=1,$ (\ref{eq8}) becomes
\begin{eqnarray*}
&&\mathbb{E}\left(\frac{\nu}{\nu+1}\right)\mathbb{E}\left(\frac{1}{\nu+1}\frac{1}{\sum_{i=0}^\nu\beta_i}\right)-
    \mathbb{E}\left(\frac{1}{\nu+1}\right)\mathbb{E}\left(\frac{\nu}{\nu+1}\frac{1}{\sum_{i=0}^\nu\beta_i}\right)\\
&&<\mathbb{E}\left(\frac{\nu}{\nu+1}\right)\mathbb{E}\left(\frac{1}{\nu+1}\right).
\end{eqnarray*}
While by Lemma \ref{lem2},\ $\sum_{i=0}^\nu\beta_i(1)\geq m_1-\frac{1}{m_1}>1,$ which implies the above inequality.

When $m_1>\lambda>1,$
\[\beta_i+(1-\lambda)\beta_i^{\prime}\leq \frac{m_1-1}{m_1-\lambda}\beta_i.\]
Combining with (\ref{eq9}), we obtain that for $1<\lambda<m_1,$
\begin{eqnarray*}
&&\mathbb{E}\left(\frac{\nu}{\nu+1}\frac{\sum_{i=0}^\nu\beta_i}{\lambda-1+\sum_{i=0}^\nu\beta_i}\right)\mathbb{E}
 \left(\frac{1}{\nu+1}\frac{\sum_{i=0}^\nu(\beta_i+(1-\lambda)\beta_i^{\prime})}{(\lambda-1+\sum_{i=0}^\nu\beta_i)^2}\right)\\
&&\ \ \ \ \leq \frac{\frac{m_1-1}{m_1-\lambda}}{m_1-\frac{\lambda}{m_1}}
             \mathbb{E}\left(\frac{\nu}{\nu+1}\frac{\sum_{i=0}^\nu\beta_i}{\lambda-1+\sum_{i=0}^\nu\beta_i}\right)
             \mathbb{E}\left(\frac{1}{\nu+1}\frac{\sum_{i=0}^\nu\beta_i}{\lambda-1+\sum_{i=0}^\nu\beta_i}\right).
\end{eqnarray*}
When $\frac{\frac{m_1-1}{m_1-\lambda}}{m_1-\frac{\lambda}{m_1}}\leq \frac{1}{\lambda}$ and $1<\lambda<m_1,$
namely $\lambda\in \left(1,\frac{m_1}{1+\sqrt{1-\frac{1}{m_1}}}\right],$ we have that
\begin{eqnarray}\label{eq10}
\nonumber &&\mathbb{E}\left(\frac{\nu}{\nu+1}\frac{\sum_{i=0}^\nu\beta_i}{\lambda-1+\sum_{i=0}^\nu\beta_i}\right)\mathbb{E}
            \left(\frac{1}{\nu+1}\frac{\sum_{i=0}^\nu(\beta_i+(1-\lambda)\beta_i^{\prime})}{(\lambda-1+\sum_{i=0}^\nu\beta_i)^2}\right)\\
&&\ \ \ \ \leq \frac{1}{\lambda}
             \mathbb{E}\left(\frac{\nu}{\nu+1}\frac{\sum_{i=0}^\nu\beta_i}{\lambda-1+\sum_{i=0}^\nu\beta_i}\right)
             \mathbb{E}\left(\frac{1}{\nu+1}\frac{\sum_{i=0}^\nu\beta_i}{\lambda-1+\sum_{i=0}^\nu\beta_i}\right).
\end{eqnarray}
Note when $1<\lambda<m_1,$
\begin{eqnarray*}
&&\sum_{i=0}^\nu(\beta_i+(1-\lambda)\beta_i^\prime)>\sum_{i=0}^\nu\beta_i>0,\\
&&\mathbb{E}\left(\frac{1}{\nu+1}\frac{\sum_{i=0}^\nu\beta_i}{\lambda-1+\sum_{i=0}^\nu\beta_i}\right)
  \mathbb{E}\left(\frac{\nu}{\nu+1}\frac{\sum_{i=0}^\nu(\beta_i+(1-\lambda)\beta_i^{\prime})}{(\lambda-1+\sum_{i=0}^\nu\beta_i)^2}\right)>0.
\end{eqnarray*}
Therefore, combining with (\ref{eq10}), we see that for $\lambda\in \left(1,\frac{m_1}{1+\sqrt{1-\frac{1}{m_1}}}\right],$
(\ref{eq8}) holds.\qed
\vskip 3mm

{\bf So far we have finished proving Theorem \ref{thm1}.}\qed

\end{document}